# A Ramanujan enigma: Part 2

Donal F. Connon

dconnon@btopenworld.com

15 February 2024

**Abstract**

We revisit Entry 17(v) in Ramanujan's Notebooks.

## 1. Introduction

A curious item appears in Ramanujan's Notebooks ([5] and [6]) which were written by him in India more than 100 years ago.

Entry 17(v) states:

If $0 < x < 1$, then

$$(1.1) \qquad \varphi(x-1) - \varphi(-x) = \pi[\gamma + \log(2\pi)]\cot(\pi x) + 2\pi \sum_{n=1}^{\infty} \sin(2n\pi x) \log n$$

where $\varphi(x)$ is defined as

$$(1.2) \qquad \varphi(x) = \sum_{n=1}^{\infty} \left[ \frac{\log n}{n} - \frac{\log(n+x)}{n+x} \right]$$

Ramanujan also notes that

$$(1.3) \qquad \sum_{n=1}^{\infty} \sin(2n\pi x) = \frac{1}{2} \cot(\pi x)$$

In his commentary on Ramanujan's Notebooks, Berndt [5, p.200] has written:

"Of course, Entry 17(v) is meaningless because the series on the right side diverges for $0 < x < 1$. In the midst of his formula, after $\cot(\pi x)$, Ramanujan inserts a parenthetical remark 'for the same limits', the meaning of which we are unable to discern." Unsurprisingly, Berndt goes on to say that (1.3) is also devoid of meaning.

We showed in [12] that Ramanujan's formulae may be rigorously expressed for $0 < x < 1$ as

$$(1.4) \qquad \varphi(x-1) - \varphi(-x) = \pi[\gamma + \log(2\pi)]\cot(\pi x) + 2\pi \lim_{s \to 0} \sum_{n=1}^{\infty} \frac{\sin(2n\pi x) \log n}{n^s}$$

and

(1.5) $$\lim_{s \to 0} \sum_{n=1}^{\infty} \frac{\sin(2n\pi x)}{n^s} = \frac{1}{2}\cot(\pi x)$$

In view of this, we now know the meaning of Ramanujan's parenthetical remark 'for the same limits'.

It was shown in [12] that

$$\varphi(x-1) - \varphi(-x) = \gamma_1(1-x) - \gamma_1(x)$$

where $\gamma_n(x)$ are the generalised Stieltjes constants defined by (1.7) below.

This results in

(1.6) $$\gamma_1(1-x) - \gamma_1(x) = \pi[\gamma + \log(2\pi)]\cot(\pi x) + 2\pi \lim_{s \to 1} \sum_{n=1}^{\infty} \frac{\sin(2n\pi x)\log n}{n^{1-s}}$$

The generalised Stieltjes constants $\gamma_n(x)$ are the coefficients in the Laurent expansion of the Hurwitz zeta function $\varsigma(s,x)$ about $s=1$

(1.7) $$\varsigma(s,x) = \sum_{n=0}^{\infty} \frac{1}{(n+x)^s} = \frac{1}{s-1} + \sum_{n=0}^{\infty} \frac{(-1)^n}{n!} \gamma_n(x)(s-1)^n$$

We have

(1.8) $$\gamma_0(x) = -\psi(x)$$

where $\psi(x)$ is the digamma function which is the logarithmic derivative of the gamma function $\psi(x) = \frac{d}{du}\log\Gamma(x)$. It is easily seen from the definition of the Hurwitz zeta function that $\varsigma(s,1) = \varsigma(s)$ and accordingly that $\gamma_n(1) = \gamma_n$.

The Stieltjes constants $\gamma_n := \gamma_n(1)$ (or the Euler-Mascheroni constants) are the coefficients of the Laurent expansion of the Riemann zeta function $\varsigma(s)$ about $s=1$.

(1.9) $$\varsigma(s) = \frac{1}{s-1} + \sum_{n=0}^{\infty} \frac{(-1)^n}{n!} \gamma_n (s-1)^n$$

Since $\lim_{s \to 1}\left[\varsigma(s) - \frac{1}{s-1}\right] = \gamma$ it is clear that $\gamma_0 = \gamma$. An elementary proof of $\gamma_0(x) = -\psi(x)$ was recently given by the author in [11] (this formula was first obtained by Berndt [4] in 1972).

This paper may usefully be read in conjunction with the author's earlier related papers, *A Ramanujan enigma involving the first Stieltjes constant* [12] and *A new representation of the Stieltjes constants* [14].

## 2. A different formulation of Ramanujan's formulae

**Theorem**

$$(2.1) \quad \lim_{s \to 0} \sum_{n=1}^{\infty} \frac{\sin(2n\pi x)}{n^s} = \frac{1}{2} \cot(\pi x)$$

**Proof**

Tričković et al. [20] showed in 2008 that for $s > 0$

$$(2.2) \quad \sum_{n=1}^{\infty} \frac{\sin(2n\pi x)}{n^s} = \frac{\pi (2\pi x)^{s-1}}{2\Gamma(s) \sin \frac{\pi s}{2}} + \sum_{n=0}^{\infty} \frac{(-1)^n \varsigma(s-2n-1)}{(2n+1)!} (2\pi x)^{2n+1}$$

and

$$(2.3) \quad \sum_{n=1}^{\infty} \frac{\cos(2n\pi x)}{n^s} = \frac{\pi (2\pi x)^{s-1}}{2\Gamma(s) \cos \frac{\pi s}{2}} + \sum_{n=0}^{\infty} \frac{(-1)^n \varsigma(s-2n)}{(2n)!} (2\pi x)^{2n}$$

It may be noted that (2.3) may be obtained from (2.2) by differentiation with respect to $x$.

Using $s\Gamma(s) = \Gamma(1+s)$ and Euler's reflection formula for the gamma function

$$\Gamma(u)\Gamma(1-u) = \frac{\pi}{\sin \pi u}$$

we may express (2.2) as

$$\sum_{n=1}^{\infty} \frac{\sin(2n\pi x)}{n^s} = \frac{(2\pi x)^{s-1} \Gamma(1+\frac{s}{2}) \Gamma(1-\frac{s}{2})}{\Gamma(1+s)} + \sum_{n=0}^{\infty} \frac{(-1)^n \varsigma(s-2n-1)}{(2n+1)!} (2\pi x)^{2n+1}$$

Hence, we obtain the limit

$$\lim_{s \to 0} \sum_{n=1}^{\infty} \frac{\sin(2n\pi x)}{n^s} = \frac{1}{2\pi x} + \sum_{n=0}^{\infty} \frac{(-1)^n \varsigma(-2n-1)}{(2n+1)!} (2\pi x)^{2n+1}$$

We have [3, p.266] in terms of the Bernoulli numbers

$$(2.4) \quad \varsigma(-m) = -\frac{B_{m+1}}{m+1}$$

and thus

$$\sum_{n=0}^{\infty}\frac{(-1)^{n}\varsigma(-2n-1)}{(2n+1)!}(2\pi x)^{2n+1}=\sum_{n=0}^{\infty}\frac{(-1)^{n+1}B_{2n+2}}{(2n+2)!}(2\pi x)^{2n+1}$$

$$=\sum_{n=1}^{\infty}\frac{(-1)^{n}B_{2n}}{(2n)!}(2\pi x)^{2n-1}$$

$$=\frac{1}{2}\sum_{n=1}^{\infty}(-1)^{n}\frac{2^{2n}B_{2n}}{(2n)!}(\pi x)^{2n-1}$$

Using the following well-known identity [15, p.130]

(2.5) $$\cot x=\frac{1}{x}+\sum_{n=1}^{\infty}(-1)^{n}\frac{2^{2n}B_{2n}}{(2n)!}x^{2n-1}\qquad,(|x|<\pi)$$

we deduce that

$$\lim_{s\to 0}\sum_{n=1}^{\infty}\frac{\sin(2n\pi x)}{n^{s}}=\frac{1}{2}\cot(\pi x)$$

**Theorem**

(2.6) $$\lim_{s\to 0}\sum_{n=1}^{\infty}\frac{\cos(2n\pi x)}{n^{s}}=-\frac{1}{2}$$

**Proof**

Letting $s\to 0$ in (2.3) gives us

$$\lim_{s\to 0}\sum_{n=1}^{\infty}\frac{\cos(2n\pi x)}{n^{s}}=\sum_{n=0}^{\infty}\frac{(-1)^{n}\varsigma(-2n)}{(2n)!}(2\pi x)^{2n}$$

and since $\varsigma(-2n)=0$ for all $n\geq 1$ we have

$$\lim_{s\to 0}\sum_{n=1}^{\infty}\frac{\cos(2n\pi x)}{n^{s}}=-\frac{1}{2}$$

since $\varsigma(0)=-\frac{1}{2}$. This interesting limit was also obtained in [12] using a different method.

**Theorem**

(2.7) $$\lim_{s\to 0}\sum_{n=1}^{\infty}(-1)^{n+1}\frac{\cos(2n\pi x)}{n^{s}}=\frac{1}{2}$$

**Proof**

Tričković et al. [20] also showed in 2008 that for $s > 0$

$$(2.8) \quad \sum_{n=1}^{\infty} (-1)^{n+1} \frac{\sin(2n\pi x)}{n^s} = \sum_{n=0}^{\infty} \frac{(-1)^n \eta(s-2n-1)}{(2n+1)!} (2\pi x)^{2n+1}$$

and

$$(2.9) \quad \sum_{n=1}^{\infty} (-1)^{n+1} \frac{\cos(2n\pi x)}{n^s} = \sum_{n=0}^{\infty} \frac{(-1)^n \eta(s-2n)}{(2n)!} (2\pi x)^{2n}$$

where $\eta(s)$ is defined by

$$(2.10) \quad \eta(s) = (1 - 2^{1-s}) \varsigma(s)$$

We have the limits

$$(2.11) \quad \lim_{s \to 0} \sum_{n=1}^{\infty} (-1)^{n+1} \frac{\sin(2n\pi x)}{n^s} = \sum_{n=0}^{\infty} \frac{(-1)^n \eta(-2n-1)}{(2n+1)!} (2\pi x)^{2n+1}$$

and

$$(2.12) \quad \lim_{s \to 0} \sum_{n=1}^{\infty} (-1)^{n+1} \frac{\cos(2n\pi x)}{n^s} = \sum_{n=0}^{\infty} \frac{(-1)^n \eta(-2n)}{(2n)!} (2\pi x)^{2n}$$

We see from (2.10) that

$$\eta(-2n) = 0 \text{ for all } n \geq 1 \quad \eta(0) = -\varsigma(0) = \frac{1}{2}$$

and thus

$$\lim_{s \to 0} \sum_{n=1}^{\infty} (-1)^{n+1} \frac{\cos(2n\pi x)}{n^s} = \frac{1}{2}$$

**Theorem**

$$(2.13) \quad \lim_{s \to 0} \sum_{n=1}^{\infty} (-1)^{n+1} \frac{\sin(n\pi x)}{n^s} = \frac{1}{2} \tan(\pi x / 2)$$

**Proof**

Using (2.8) we see that

$$\lim_{s \to 0} \sum_{n=1}^{\infty} (-1)^{n+1} \frac{\sin(2n\pi x)}{n^s} = \sum_{n=0}^{\infty} \frac{(-1)^n (1-2^{2n+2}) \varsigma(-2n-1)}{(2n+1)!} (2\pi x)^{2n+1}$$

$$= \sum_{n=0}^{\infty} \frac{(-1)^{n+1} (1-2^{2n+2}) B_{2n+2}}{(2n+2)!} (2\pi x)^{2n+1}$$

$$= \sum_{m=1}^{\infty} \frac{(-1)^{m+1} (2^{2m} - 1) B_{2m}}{(2m)!} (2\pi x)^{2m-1}$$

Using (2.5) and $\tan x = \cot x - 2\cot(2x)$ we obtain the well-known result

$$\tan x = \sum_{m=1}^{\infty} \frac{(-1)^{m+1}(2^{2m}-1)2^{2m} B_{2m}}{(2m)!} x^{2m-1}$$

and hence we have

$$\lim_{s \to 0} \sum_{n=1}^{\infty} (-1)^{n+1} \frac{\sin(n\pi x)}{n^s} = \frac{1}{2}\tan(\pi x/2)$$

**Theorem**

(2.14) $$B_{2m+1}(x) = -\frac{(2m+1)}{2} x^{2m} + \sum_{n=0}^{m} \binom{2m+1}{2n+1} B_{2(m-n)} x^{2n+1}$$

**Proof**

As noted by Tričković et al. [20] letting $s = 2m+1$ in (2.2) results in

$$\sum_{n=1}^{\infty} \frac{\sin(2n\pi x)}{n^{2m+1}} = \frac{(-1)^m \pi (2\pi x)^{2m}}{2(2m)!} + \sum_{n=0}^{m} \frac{(-1)^n \varsigma(2m-2n)}{(2n+1)!} (2\pi x)^{2n+1}$$

where the series truncates because $\varsigma(-2k) = 0$ for all $k \geq 1$. Using Euler's formula

$$\varsigma(2k) = (-1)^{k+1} \frac{B_{2k}(2\pi)^{2k}}{2(2k)!}$$

we see that

$$\sum_{n=0}^{m} \frac{(-1)^n \varsigma(2m-2n)}{(2n+1)!} (2\pi x)^{2n+1} = \frac{1}{2}(2\pi)^{2m+1} \sum_{n=0}^{m} \frac{(-1)^{m+1}}{(2n+1)!} \frac{B_{2(m-n)}}{(2m-2n)!} x^{2n+1}$$

Since it is known that [2]

$$\sum_{n=1}^{\infty} \frac{\sin(2n\pi x)}{n^{2m+1}} = (-1)^{m+1} \frac{(2\pi)^{2m+1}}{2(2m+1)!} B_{2m+1}(x)$$

we have

$$(-1)^{m+1} \frac{(2\pi)^{2m+1}}{2(2m+1)!} B_{2m+1}(x) = \frac{(-1)^m \pi (2\pi x)^{2m}}{2(2m)!} + \frac{1}{2}(2\pi)^{2m+1} \sum_{n=0}^{m} \frac{(-1)^{m+1}}{(2n+1)!} \frac{B_{2(m-n)}}{(2m-2n)!} x^{2n+1}$$

which simplifies to

$$B_{2m+1}(x) = -\frac{(2m+1)}{2}x^{2m} + \sum_{n=0}^{m}\binom{2m+1}{2n+1}B_{2(m-n)}\,x^{2n+1}$$

Since $B'_n(x) = nB_{n-1}(x)$ differentiation gives us

$$B_{2m}(x) = -mx^{2m-1} + \frac{1}{2m+1}\sum_{n=0}^{m}\binom{2m+1}{2n+1}(2n+1)B_{2(m-n)}\,x^{2n}$$

$$= -mx^{2m-1} + \sum_{n=0}^{m}\binom{2m}{2n}B_{2n}\,x^{2m-2n}$$

which appears in [sr, p.60]. These formulae are clearly related to the more familiar identity [17, p.59]

$$B_m(x) = \sum_{n=0}^{m}\binom{m}{n}B_n\,x^{m-n}$$

Letting $m = 0$ gives us

$$\sum_{n=1}^{\infty}\frac{\sin(2n\pi x)}{n} = \frac{\pi}{2} + 2\pi x\varsigma(0)$$

and we obtain the well-known Fourier Series

$$\sum_{n=1}^{\infty}\frac{\sin(2n\pi x)}{n} = \pi(\tfrac{1}{2} - x)$$

**Theorem**

(2.15) $\quad \gamma_1(1-x) - \gamma_1(x) = -\dfrac{\gamma + \log(2\pi x)}{x} + \pi[\gamma + \log(2\pi)]\cot(\pi x)$

$$+2\pi\sum_{n=0}^{\infty}\frac{(-1)^{n+1}\varsigma'(-2n-1)}{(2n+1)!}(2\pi x)^{2n+1}$$

where $\gamma_1(x)$ is the first Stieltjes constant.

**Proof**

We recall (2.2) in the format

(2.16) $\quad \displaystyle\sum_{n=1}^{\infty}\frac{\sin(2n\pi x)}{n^s} = \frac{(2\pi x)^{s-1}\Gamma(1+\frac{s}{2})\Gamma(1-\frac{s}{2})}{\Gamma(1+s)} + \sum_{n=0}^{\infty}\frac{(-1)^n\varsigma(s-2n-1)}{(2n+1)!}(2\pi x)^{2n+1}$

Defining $f(s)$ by

$$f(s) = \frac{(2\pi x)^{s-1}\Gamma(1+\tfrac{s}{2})\Gamma(1-\tfrac{s}{2})}{\Gamma(1+s)}$$

we have the first derivative

$$f'(s) = [\log(2\pi x) + \tfrac{1}{2}\psi(1+\tfrac{s}{2}) - \tfrac{1}{2}\psi(1-\tfrac{s}{2}) - \psi(1+s)]f(s)$$

Higher derivatives of $f(s)$ may therefore be compactly expressed in terms of the (exponential) complete Bell polynomials [9].

Differentiation of (2.16) with respect to $s$ results in

$$-\sum_{n=1}^{\infty} \frac{\log n}{n^s} \sin(2n\pi x) = [\log(2\pi x) + \tfrac{1}{2}\psi(1+\tfrac{s}{2}) - \tfrac{1}{2}\psi(1-\tfrac{s}{2}) - \psi(1+s)]f(s)$$

$$+ \sum_{n=0}^{\infty} \frac{(-1)^n \varsigma'(s-2n-1)}{(2n+1)!}(2\pi x)^{2n+1}$$

and the limit $s \to 0$ gives us

$$\lim_{s \to 0}\sum_{n=1}^{\infty} \frac{\log n}{n^s}\sin(2n\pi x) = -\frac{\gamma + \log(2\pi x)}{2\pi x} + \sum_{n=0}^{\infty}\frac{(-1)^{n+1}\varsigma'(-2n-1)}{(2n+1)!}(2\pi x)^{2n+1}$$

This is in agreement with equation (9.12) in [14] which was obtained in a completely different manner.

We also showed in [12] that

(2.17) $\qquad \gamma_1(1-x) - \gamma_1(x) = 2\pi \lim_{s \to 0}\sum_{n=1}^{\infty} \frac{\log n}{n^s}\sin(2n\pi x) + \pi[\gamma + \log(2\pi)]\cot(\pi x)$

We therefore obtain

$$\gamma_1(1-x) - \gamma_1(x) = -\frac{\gamma + \log(2\pi x)}{x} + \pi[\gamma + \log(2\pi)]\cot(\pi x)$$

$$+ 2\pi\sum_{n=0}^{\infty}\frac{(-1)^{n+1}\varsigma'(-2n-1)}{(2n+1)!}(2\pi x)^{2n+1}$$

$$= -\frac{\gamma + \log(2\pi x)}{x} + \pi[\gamma + \log(2\pi)]\cot(\pi x)$$

$$+ 2\pi\sum_{k=1}^{\infty}\frac{(-1)^k \varsigma'(1-2k)}{(2k-1)!}(2\pi x)^{2k-1}$$

One could also substitute

$$2k\varsigma'(1-2k) = \frac{\varsigma'(2k)}{\varsigma(2k)} B_{2k} + \left[H^{(1)}_{2k-1} - \gamma - \log 2\pi\right] B_{2k}$$

Letting $x = \frac{1}{2}$ gives us

(2.18) $$\sum_{n=0}^{\infty} \frac{(-1)^{n+1} \varsigma'(-2n-1)}{(2n+1)!} (\pi)^{2n+1} = \frac{\gamma + \log \pi}{\pi}$$

**Theorem**

(2.19) $$\sum_{n=1}^{\infty} \frac{\log n}{n} \cos 2\pi nt = \frac{1}{2}\left[\varsigma''(0,t) + \varsigma''(0,1-t)\right] + [\gamma + \log(2\pi)]\log(2\sin \pi t)$$

**Proof**

Let us now <u>assume</u> that we may validly integrate (2.17) through the limit, term by term over the interval $[t, \frac{1}{2}]$.

We employ the integral [12]

$$\int_1^t \gamma_n(x) = \frac{(-1)^{n+1}}{n+1} [\varsigma^{(n+1)}(0,t) - \varsigma^{(n+1)}(0)]$$

and thus

$$\int_{\frac{1}{2}}^{t} \gamma_1(x)dx = \int_1^t \gamma_1(x)dx - \int_1^{\frac{1}{2}} \gamma_1(x)dx$$

$$= \frac{1}{2}[\varsigma''(0,t) - \varsigma''(0,\tfrac{1}{2})]$$

We have

$$\int_{\frac{1}{2}}^{t} \gamma_1(1-x)dx = -\int_{\frac{1}{2}}^{1-t} \gamma_1(u)du$$

$$= -\frac{1}{2}[\varsigma''(0,1-t) - \varsigma''(0,\tfrac{1}{2})]$$

This gives us

$$\sum_{n=1}^{\infty} \frac{\log n}{n} \cos 2\pi nt = \sum_{n=1}^{\infty} (-1)^n \frac{\log n}{n} + [\gamma + \log(2\pi)]\log(\sin \pi t)$$

$$+ \frac{1}{2}[\varsigma''(0,t) + \varsigma''(0,1-t)] - \varsigma''(0,\tfrac{1}{2})$$

It is well known that

$$\sum_{n=1}^{\infty}(-1)^{n+1}\frac{\log n}{n}=\frac{1}{2}\log^{2}2-\gamma\log 2$$

Using the well-known formula $\varsigma(s,\tfrac{1}{2})=(2^{s}-1)\varsigma(s)$ we readily find that

$$\varsigma''\left(0,\tfrac{1}{2}\right)=2\varsigma'(0)\log 2-\frac{1}{2}\log^{2}2$$

and Lerch's formula

$$\varsigma'(0,x)=\log\Gamma(x)-\frac{1}{2}\log(2\pi)$$

gives us

$$\varsigma''\left(0,\tfrac{1}{2}\right)=-\log(2\pi)\log 2-\frac{1}{2}\log^{2}2$$

Hence we obtain

$$\sum_{n=1}^{\infty}\frac{\log n}{n}\cos 2\pi nt=\frac{1}{2}\left[\varsigma''(0,t)+\varsigma''(0,1-t)\right]+\left[\gamma+\log(2\pi)\right]\log(2\sin\pi t)$$

This identity was derived by Ramanujan and it appears, albeit in a highly disguised form, in Berndt's book [5, Part I, p.208] and [7b]. It was reported by Deninger [16] in 1984 and was also derived by Blagouchine [7] in 2014.

We showed in [14] that

(2.20) $$\varsigma''(0,t)=\sum_{n=1}^{\infty}\frac{[\gamma+\log(2n\pi)]^{2}\sin(2n\pi t)}{n\pi}-\frac{1}{4}\varsigma(2)\sum_{n=1}^{\infty}\frac{\sin(2n\pi t)}{n\pi}$$

$$+\sum_{n=1}^{\infty}\frac{[\gamma+\log(2\pi n)]\cos(2n\pi t)}{n}$$

**Theorem**

(2.21) $$\lim_{s\to 0}2\sum_{n=1}^{\infty}\frac{\log n}{n^{s}}\cos(2n\pi x)=\psi(x)+\frac{\pi}{2}\cot(\pi x)+\gamma+\log(2\pi)$$

**Proof**

We express (2.2) in the form

(2.22) $$\sum_{n=1}^{\infty}\frac{\cos(2n\pi x)}{n^{s}}=\frac{(2\pi x)^{s-1}\Gamma(\tfrac{1+s}{2})\Gamma(\tfrac{1-s}{2})}{2\Gamma(s)}+\sum_{n=0}^{\infty}\frac{(-1)^{n}\varsigma(s-2n)}{(2n)!}(2\pi x)^{2n}$$

Letting

$$g(s) = \frac{(2\pi x)^{s-1}\Gamma(\frac{1+s}{2})\Gamma(\frac{1-s}{2})}{2\Gamma(s)}$$

we have the first derivative

$$g'(s) = [\log(2\pi x) + \tfrac{1}{2}\psi(\tfrac{1+s}{2}) - \tfrac{1}{2}\psi(\tfrac{1-s}{2}) - \psi(s)]g(s)$$

As before, higher derivatives of $g(s)$ may be compactly expressed in terms of the (exponential) complete Bell polynomials [9].

Differentiating (2.22) results in

$$-\sum_{n=1}^{\infty}\frac{\log n}{n^s}\cos(2n\pi x) = [\log(2\pi x) + \tfrac{1}{2}\psi(\tfrac{1+s}{2}) - \tfrac{1}{2}\psi(\tfrac{1-s}{2}) - \psi(s)]\frac{(2\pi x)^{s-1}\Gamma(\frac{1+s}{2})\Gamma(\frac{1-s}{2})}{2\Gamma(s)}$$

$$+\sum_{n=0}^{\infty}\frac{(-1)^n\varsigma'(s-2n)}{(2n)!}(2\pi x)^{2n}$$

and we have the limit

$$-\lim_{s\to 0}\sum_{n=1}^{\infty}\frac{\log n}{n^s}\cos(2n\pi x) = -\lim_{s\to 0}\psi(s)\frac{(2\pi x)^{s-1}\Gamma(\frac{1+s}{2})\Gamma(\frac{1-s}{2})}{2\Gamma(s)}$$

$$+\sum_{n=0}^{\infty}\frac{(-1)^n\varsigma'(-2n)}{(2n)!}(2\pi x)^{2n}$$

Using $\Gamma(\tfrac{1}{2}) = \sqrt{\pi}$ and the limit

$$\lim_{s\to 0}\frac{\psi(s)}{\Gamma(s)} = \lim_{s\to 0}\frac{s\psi(s)}{\Gamma(1+s)}$$

$$= \lim_{s\to 0}\frac{s[\psi(1+s) - \tfrac{1}{s}]}{\Gamma(1+s)} = -1$$

we obtain

$$\lim_{s\to 0}\sum_{n=1}^{\infty}\frac{\log n}{n^s}\cos(2n\pi x) = -\frac{1}{4x} - \sum_{n=0}^{\infty}\frac{(-1)^n\varsigma'(-2n)}{(2n)!}(2\pi x)^{2n}$$

It is well-known that [17, p.9] for $n \geq 1$

$$\varsigma'(-2n) = (-1)^n\frac{(2n)!}{2(2\pi)^{2n}}\varsigma(2n+1)$$

and we have [17]

$$\varsigma'(0) = -\frac{1}{2}\log(2\pi)$$

Thus, we obtain

(2.23) $$\lim_{s \to 0} 2\sum_{n=1}^{\infty} \frac{\log n}{n^s} \cos(2n\pi x) = -\frac{1}{2x} + \log(2\pi) - \sum_{n=1}^{\infty} \varsigma(2n+1) x^{2n}$$

We have [17, p.160] for $|x| < |a|$

$$\sum_{n=1}^{\infty} \varsigma(2n+1, a) x^{2n} = -\frac{1}{2}[\psi(a+x) + \psi(a-x)] + \psi(a)$$

and with $a = 1$ this becomes

$$\sum_{n=1}^{\infty} \varsigma(2n+1) x^{2n} = -\frac{1}{2}[\psi(1+x) + \psi(1-x)] - \gamma$$

Employing [17, p.14] for $x > 0$

$$\psi(1+x) - \psi(1-x) = \frac{1}{x} - \pi \cot \pi x$$

we deduce that

(2.24) $$\lim_{s \to 0} 2\sum_{n=1}^{\infty} \frac{\log n}{n^s} \cos(2n\pi x) = \psi(x) + \frac{\pi}{2} \cot(\pi x) + \gamma + \log(2\pi)$$

which was previously shown in [14]. This may also be expressed as

(2.25) $$\lim_{s \to 0} 2\sum_{n=1}^{\infty} \frac{\log n}{n^s} \cos(2n\pi x) = \gamma + \log(2\pi) + \frac{1}{2}[\psi(x) - \psi(1-x)]$$

Let us now assume that we may validly integrate (2.25) through the limit, term by term over the interval $[t, \frac{1}{2}]$. This gives us

(2.26) $$\sum_{n=1}^{\infty} \frac{\log n}{n} \sin 2\pi nt = \frac{\pi}{2}\left[\log \Gamma(t) - \log \Gamma(1-t)\right] + [\gamma + \log(2\pi)]\pi\left(t - \frac{1}{2}\right)$$

which was previously reported by Deninger [16]. It is easily seen that this is equivalent to Kummer's Fourier series for $\log \Gamma(t)$ [4].

**Theorem**

(2.1) $$\sum_{n=1}^{\infty}\frac{\sin(2n\pi x)}{n^{2m}}=(-1)^m\frac{(2\pi x)^{2m-1}}{(2m-1)!}[\gamma+\log(2\pi x)-\psi(2m)]$$

$$+\sum_{n=0}^{m-2}\frac{(-1)^n\varsigma(2m-2n-1)}{(2n+1)!}(2\pi x)^{2n+1}+\sum_{n=m}^{\infty}\frac{(-1)^n\varsigma(2m-2n-1)}{(2n+1)!}(2\pi x)^{2n+1}$$

**Proof**

We let $s\to 2m$ in (2.2) to obtain

(2.27) $$\sum_{n=1}^{\infty}\frac{\sin(2n\pi x)}{n^{2m}}=\lim_{s\to 2m}\left[\frac{\pi(2\pi x)^{s-1}}{2\Gamma(s)\sin\frac{\pi s}{2}}+\frac{(-1)^{m-1}\varsigma(s-2m+1)}{(2m-1)!}(2\pi x)^{2m-1}\right]$$

$$+\sum_{n=0}^{m-2}\frac{(-1)^n\varsigma(2m-2n-1)}{(2n+1)!}(2\pi x)^{2n+1}+\sum_{n=m}^{\infty}\frac{(-1)^n\varsigma(2m-2n-1)}{(2n+1)!}(2\pi x)^{2n+1}$$

where we have separated out the term corresponding to $n=m-1$ and started the infinite series at $n=m$.

We write the term

$$\frac{\pi(2\pi x)^{s-1}}{2\Gamma(s)\sin\frac{\pi s}{2}}+\frac{(-1)^{m-1}\varsigma(s-2m+1)}{(2m-1)!}(2\pi x)^{2m-1}$$

$$=\left[\frac{\pi(2\pi x)^{s-1}}{2\Gamma(s)}+(-1)^{m-1}\frac{\sin\frac{\pi s}{2}}{s-2m}\frac{(s-2m)\varsigma(s-2m+1)}{(2m-1)!}(2\pi x)^{2m-1}\right]/\sin\frac{\pi s}{2} \quad (*)$$

$$:=L(s)/\sin\frac{\pi s}{2}$$

Noting that

$$\lim_{s\to 2m}\frac{\sin\frac{\pi s}{2}}{s-2m}=\lim_{s\to 2m}\frac{\pi}{2}\cos\frac{\pi s}{2}=(-1)^m\frac{\pi}{2}$$

and

$$\lim_{s\to 2m}(s-2m)\varsigma(s-2m+1)=\lim_{s\to 1}(s-1)\varsigma(s)=1$$

we deduce that

$$\lim_{s\to 2m}\left[\frac{\pi(2\pi x)^{s-1}}{2\Gamma(s)}+(-1)^{m-1}\frac{\sin\frac{\pi s}{2}}{s-2m}\frac{(s-2m)\varsigma(s-2m+1)}{(2m-1)!}(2\pi x)^{2m-1}\right]=0$$

Hence, we may apply L'Hôpital's rule to (*) and this gives us

$$L'(s) = \frac{\pi(2\pi x)^{s-1}}{2}\frac{\log(2\pi x) - \psi(s)}{\Gamma(s)}$$

$$+(-1)^{m-1}\frac{(2\pi x)^{2m-1}}{(2m-1)!}\left[\phi(s)\frac{d}{ds}[(s-2m)\varsigma(s-2m+1)] + \phi'(s)(s-2m)\varsigma(s-2m+1)\right]$$

where $\phi(s) := \frac{\sin\frac{\pi s}{2}}{s-2m}$. We express $\sin\frac{\pi s}{2} = \sum_{n=1}^{\infty}\alpha_n(s-2m)^n$ and easily deduce that

$$\alpha_1 = (-1)^m\frac{\pi}{2} \text{ and } \alpha_2 = 0$$

Hence

$$L'(2m) = \frac{\pi(2\pi x)^{2m-1}}{2}\frac{\log(2\pi x) - \psi(2m)}{\Gamma(2m)} - \frac{\pi}{2}\frac{(2\pi x)^{2m-1}}{(2m-1)!}\frac{d}{ds}(s-2m)\varsigma(s-2m+1)\bigg|_{s=2m}$$

We see from () that

$$\frac{d}{ds}(s-2m)\varsigma(s-2m+1)\bigg|_{s=2m} = -\gamma$$

we have

$$L'(2m) = \frac{\pi}{2}\frac{(2\pi x)^{2m-1}}{(2m-1)!}[\gamma + \log(2\pi x) - \psi(2m)]$$

Hence, we obtain

$$\sum_{n=1}^{\infty}\frac{\sin(2n\pi x)}{n^{2m}} = (-1)^m\frac{(2\pi x)^{2m-1}}{(2m-1)!}[\gamma + \log(2\pi x) - \psi(2m)]$$

$$+ \sum_{n=0}^{m-2}\frac{(-1)^n\varsigma(2m-2n-1)}{(2n+1)!}(2\pi x)^{2n+1} + \sum_{n=m}^{\infty}\frac{(-1)^n\varsigma(2m-2n-1)}{(2n+1)!}(2\pi x)^{2n+1}$$

The above is a simplified derivation of the proof provided by Tričković et al. [18].

Via a lengthy proof, Tričković et al. [18] used (2.2) and (2.3) to show that

$$\frac{(-1)^m(2m-1)!}{(2\pi)^{2m-1}}\sum_{n=1}^{\infty}\frac{\sin(2n\pi x)}{n^{2m}} = x^{2m-1}\log x + \varsigma'(1-2m,1-x) - \varsigma'(1-2m,1+x)$$

$$\frac{(-1)^m(2m-2)!}{(2\pi)^{2m-2}}\sum_{n=1}^{\infty}\frac{\cos(2n\pi x)}{n^{2m-1}} = x^{2m-2}\log x - \varsigma'(2-2m,1-x) - \varsigma'(2-2m,1+x)$$

These results were previously obtained by Adamchik [1] for $0 < x < 1$ in the form

$$\varsigma'(-m, x) + (-1)^m \varsigma'(-m, 1-x) = \pi i \frac{B_{m+1}}{m+1} + e^{-i\pi m/2} \frac{m!}{(2\pi)^m} Li_{m+1}(e^{2\pi i x})$$

and their equivalence may be seen by noting that

$$\varsigma(s, 1+x) = \varsigma(s, x) - \frac{1}{x^s}$$

so that

$$\varsigma'(1-2m, 1+x) = \varsigma'(1-2m, x) + x^{2m-1} \log x$$

## 3. Acknowledgement

I thank Roy J Hughes for forwarding me a copy of [18].

## 4. Open access to our own work

This paper contains references to various other papers and, rather surprisingly, most of them are currently freely available on the internet. Surely now is the time that <u>all</u> of <u>our</u> work should be freely accessible by <u>all</u>. The mathematics community should lead the way on this by publishing <u>everything</u> on arXiv, or in an equivalent open access repository. We think it, we write it, so why hide it? You know it makes sense.

## REFERENCES

[1]   V.S. Adamchik, A Class of Logarithmic Integrals. Proceedings of the 1997
      International Symposium on Symbolic and Algebraic Computation.
      ACM, Academic Press, 1-8, 2001.
      https://www.researchgate.net/publication/221564534_A_Class_of_Logarithmic_Integrals

[2]   T.M. Apostol, Mathematical Analysis, Second Ed., Addison-Wesley Publishing
      Company, Menlo Park (California), London and Don Mills (Ontario), 1974.

[3]   T.M. Apostol, Introduction to Analytic Number Theory.
      Springer-Verlag, New York, Heidelberg and Berlin, 1976.

[4]   B.C. Berndt, On the Hurwitz Zeta-function.
      Rocky Mountain Journal of Mathematics, vol. 2, no. 1, pp. 151-157, 1972.

[5]   B.C. Berndt, Ramanujan's Notebooks. Part I, Springer-Verlag, 1989.

[6]   B.C. Berndt, Chapter eight of Ramanujan's Second Notebook.
      J. Reine Agnew. Math, Vol. 338, 1-55, 1983.
      http://www.digizeitschriften.de/no_cache/home/open-access/nach-zeitschriftentiteln/

[7] I. V. Blagouchine, Rediscovery of Malmstén's integrals, their evaluation by contour integration methods and some related results.
The Ramanujan Journal, vol. 35, no. 1, pp. 21–110 (2014).
https://www.researchgate.net/profile/Iaroslav_Blagouchine/publications

[8] I.V. Blagouchine, A theorem for the closed-form evaluation of the first generalized Stieltjes constant at rational arguments and some related summations.
Journal of Number Theory 148 (2015) 537-592.
http://dx.doi.org/10.1016/j.jnt.2014.08.009
http://arxiv.org/pdf/1401.3724.pdf

[9] D.F. Connon, Various applications of the (exponential) complete Bell polynomials. 2010.
https://arxiv.org/abs/1001.2835

[10] D.F. Connon, A formula connecting the Bernoulli numbers with the Stieltjes constants. 2011.
https://arxiv.org/abs/1104.4772

[11] D.F. Connon, Some integrals and series involving the Stieltjes constants, 2018.
https://arxiv.org/abs/1801.05711

[12] D.F. Connon, A Ramanujan enigma involving the first Stieltjes constant, 2019.
https://arxiv.org/abs/1901.03382

[13] D.F. Connon, Some new formulae involving the Stieltjes constants, 2019.
https://arxiv.org/abs/1902.00510

[14] D.F. Connon, A new representation of the Stieltjes constants, 2022.
http://arxiv.org/abs/2201.05084

[15] G. Boros and V.H. Moll, Irresistible Integrals: Symbolics, Analysis and Experiments in the Evaluation of Integrals. Cambridge University Press, 2004.

[16] C. Deninger, On the analogue of the formula of Chowla and Selberg for real quadratic fields. J. Reine Angew. Math., 351 (1984), 172–191.
http://www.digizeitschriften.de/dms/toc/?PPN=PPN243919689_0351

[17] H.M. Srivastava and J. Choi, Series Associated with the Zeta and Related Functions. Kluwer Academic Publishers, Dordrecht, the Netherlands, 2001.

[18] S.B. Tričković & M.S. Stanković On the closed form of Clausen functions, Integral Transforms and Special Functions, 2023, 34:6, 469-477.
https://www.researchgate.net/profile/Slobodan-Trickovic

[19] S.B. Tričković, M.S. Stanković, M.V. Vidanović, On the summation of Schlömilch series. Int. Trans. Spec. Func. 2020, 31(5), 339–367.
https://www.researchgate.net/profile/Slobodan-Trickovic

[20] S.B. Tričković, M.V. Vidanović, M.S. Stanković, On the Summation of Trigonometric Series. Integral Transforms and Special Functions, 2008, 19:6, 441-452.
https://www.researchgate.net/profile/Slobodan-Trickovic

Wessex House,
Devizes Road,
Upavon,
Wiltshire SN9 6DL